\documentclass[12pt]{article}
\usepackage{amsmath}
\usepackage{amssymb}
\usepackage{latexsym}
\usepackage{epsfig}
\usepackage{graphicx}
\usepackage{mathptmx}
\usepackage[mtpcal]{}
\usepackage{mathtools}
\usepackage{physics}
\usepackage{cancel}
\usepackage{subcaption}
\begin{document}
\setlength{\baselineskip}{7mm}

\begin{center}
{\huge Centrality Estimators for Probability Density Functions}         % used by \maketitle
\end{center}
%\vspace*{5mm}
\begin{center}
\begin{tabular}{c}
Djemel Ziou \\
D\'epartement d'informatique  \\
 Universit\'e de Sherbrooke  \\
 Sherbrooke, Qc., Canada J1K 2R1\\
 Djemel.Ziou@usherbrooke.ca
%Formulation made in 2017 and revised in 2022.
\end{tabular}
\end{center}

\setlength{\baselineskip}{7mm}

%\vspace*{20mm}
%{\bf Corresponding author:} Djemel.Ziou@usherbrooke.ca

\vspace*{5mm}

\begin{abstract}
In this report,  we explore the data selection leading to  a family  of estimators maximizing a centrality. The family allows a nice properties leading to accurate and robust probability density function fitting according to some criteria we define. We  establish a link between the centrality estimator and the maximum likelihood,  showing that the latter is a particular case. Therefore, a new probability interpretation of Fisher maximum likelihood is provided. We will introduce and study two specific centralities that we have named H\"older and Lehmer estimators. A numerical simulation is provided showing the effectiveness of the proposed families of estimators opening the door to development of new concepts and algorithms in  machine learning, data mining, statistics, and data analysis.
\end{abstract}


%\end{frontmatter}

%\linenumbers

%% main text
\section{Introduction}
Maximum likelihood estimator (MLE) is a general statistical estimation method proposed by R. A. Fisher almost one century before. Since then,  many writings have been published about the MLE, its variants, its properties, its interpretation,  and  the controversy it provokes. 
 In the case of iid data, Fisher defined the probability of $X=\{x_1, \cdots, x_n\}$ occurring as being proportional to the probability $P_\theta(X) = \prod_i^n P_\theta(x_i)$. For the discrete random variable $x$,    the probability $P_\theta (x_i)=P_\theta(x=x_i)$. In the continuous case, assuming  that the iid data  $X$ are sampled from a probability density function (PDF) $h(x|\theta)$ with a parameter $\theta$,  a value $x_i$ is assumed to be known within some accuracy of measurement $\epsilon$. The probability  that $y$ falls in the cell $[x_i - \epsilon/2, x_i+\epsilon/2]$ is $P_\theta(x_i)=\int_{y \in [x_i - \epsilon/2, x_i+\epsilon/2]} h(y|\theta) dy$. Let assume that the cell width  $\epsilon$ is small and the PDF  $h(y|\theta)$ is almost constant within the interval $[x_i - \epsilon/2, x_i+\epsilon/2]$, it follows that $P_\theta(x_i) \simeq  \epsilon h(x_i|\theta)$.  Whatever $x$ is discrete or continuous random variable,  the most probable set of values for $\theta$ will make  the probability $P_\theta(X)$ maximal. 
 
 Let us take things from a different perspective. Our motivation is caused on the one hand, by errors in the real data, the iid hypothesis, the unknown PDF, the size of the cell $\epsilon$. The errors in the data can have different origins including the measuring instruments, the human, the unavailability of the measurements. It is  often difficult to formally sustain  that the data are assumed to be generated from a known PDF or they are iid. In practice, the iid assumption is at the basis of the  likelihood definition and  the chosen PDF is retained  if it leads to acceptable clustering, segmentation, prediction or decision.    In the case of a continuous random variable, the cell size $\epsilon$ is considered as a constant, with no effect on the maximum likelihood estimate. The transition from the probability to a PDF is valid if the cell size is chosen adequately. The information theory and statistical model selection made it clear that the cell size  depends on  the observed Fisher information~\cite{Ziou09}. In addition to validity of transition from probability to PDF, the cell size $\epsilon$ is of great importance because it informs about the amount of required data for the PDF fitting. 
 
%which may have different sources including  badly calibrated  measuring  instrument, produced by a human, or %retrieved from available a data collection. 
%\begin{itemize}
%\item  badly calibrated instruments, limited resources to acquire sufficient data, lost of a part of the data,  misunderstanding of the measured phenomenon, and even  bugs in computer program leads to produce corrupted data.
%\item the PDF from which the data are generated is often unknown. 
%it is  often difficult to formally sustain  that the data are assumed to be generated from a known PDF. In %practice, the chosen PDF is retained  if it leads to acceptable clustering, segmentation, prediction or %decision. 
%\item IID assumption
%\item 
% EPSILON DEPEND DE THETA 
%\end{itemize}

%This reality leads to a reconsideration of likelihood. Indeed, the most probable values of the $\theta$ which maximize $P_\theta()$ must be robust with respect to the bias that may be introduced by this reality. 

On the other side,  the likelihood as the product of probabilities $\prod_i^n P_\theta(x_i)$ is somehow "harsh" as it may be highly corrupted by a single low outlier. Considering that the root $n^{th}$ of this product is a geometric mean indicates that  $\prod_i^n P_\theta(x_i)$ is closer to the smallest values of  $P_\theta(x_i)$. However,
accepting that the most probable  values for the $\theta$'s is observed on any cell over the support covered by  $X$ leads to consider $P_\theta(X)$ as a centrality measure of $P_\theta(x_i)$.  Among the many existing centralities, we choose those of H\"older and Lehmer, given respectively by:
\begin{equation}
H_\theta(X; \alpha)=(\sum_{i=1}^n  \lambda_i P_\theta^\alpha(x_i))^{1/\alpha}
\label{Hmean}
\end{equation}
\begin{equation}
L_\theta(X,\alpha)=\frac{\sum_{i=1}^n  \lambda_i P_\theta^\alpha(x_i)} {\sum_{i=1}^n  \lambda_iP_\theta^{\alpha-1}(x_i)}
\label{Lmean}
\end{equation}
where the non-negative weights fulfills  $\sum_{i=1}^n  \lambda_i = 1$.
Each of H\"older and Lehmer centralities implement specific data selection criteria~\cite{Ziou23b}. Note that the uses in data selection, estimation and modeling in various fields of these  centralities and Pythagorean means which are special cases are given in~\cite{Ziou23a,Ziou23b}. 
The centralities  has nice properties allowing to overcome the above mentioned issue by choosing the $\alpha$ leading to a higher probability $P_\theta(X)$.  The estimation can be more accurate, according to some criteria,  at the cost of  increased complexity.  From now on,  as in the likelihood case, we are not interested on the probability of the data $P_\theta(X)$, rather our focus is the parameter $\theta$. This generalization of the likelihood  have several advantages: 1) it overcomes the IID assumption; 2) thanks to embedded data selection criterion in each of the centralities, the  maximum of a centrality  (which we call centrality estimator and we denote C-estimator), if exists, can be less sensitive to  errors in data. However, this generalization leads to new issues and a better understanding of  Fisher's proposal. The properties of $P_\alpha(X)$, the existence of the maximum, Fisher information, and the selection of the suitable centrality  are among issues. 

Our goal is not to provide a rigorous mathematical reasoning, but rather to deal with the computation of C-estimators, their properties and their evaluation, under the following assumptions: 1)  the random variable is  considered continuous and the values of the PDF $h(x_1|\theta), \cdots, h(x_n|\theta)$ are not all equal as it is the case of  the uniform distribution; 2) The centralities in Eqs.~\ref{Hmean} and~\ref{Lmean} are positive, continuous and derivable w.r.t to both $\alpha$ and $\theta$. The next two sections are dedicated respectively to the study of H\"older and Lehmer centralities.
The sections~\ref{HMeasureS} and~\ref{LMeasureS}  describe H\"older and Lehmer centrality measures and some of their properties. A comparison of the two measures is described in section~\ref{CHL}. Section~\ref{MCE} is devoted to the maximum centrality estimation. 
The accuracy of the C-estimators is explained in section~\ref{AE}. A case study is the subject of the section~\ref{CS}.

\section{H\"older measure}
\label{HMeasureS}
In this work, we consider the continuous random variable case where the cell size is function of the Fisher information as it is explained before; $P_\theta(y_i) \simeq  \epsilon(\theta) h(y_i|\theta)$. In this case, the centrality, denoted H-Centrality (H-C)   when using Eq.~\ref{Hmean}  is rewritten  as follow: 
\begin{equation}
 H_\alpha (\theta) =  (\sum_i \lambda_i ( \epsilon(\theta) h(x_i|\theta))^\alpha)^{1/ \alpha}
\label{HCentrality}
\end{equation}
where $\alpha \in \mathbb{R}$. We assume that  $ H_\alpha (\theta)>0$ is continuous and derivable w.r.t to both $\alpha$ and $\theta \in \Theta$ and $ \epsilon (\theta)>0$ is continuous and derivable w.r.t to $\theta$.  An interpretation of this equation is as follows: a data element $h(x_i|\theta)$ is transformed to $h^\alpha(x_i|\theta)$, then the arithmetic mean is computed in the transformed space and  it is back transformed to the original space.  As it will be explained further, the transform and its inverse implement data selection criteria.  We will now state properties of $H_\alpha(\theta)$.   

{\bf P1:} $lim_{\alpha \to 0}  H_\alpha(\theta)=\prod_i (\epsilon(\theta) h(x_i|\theta))^{\lambda_i}$. The limit of $\lim_{\alpha \to 0} ln  H_\alpha(\theta) = 0/0$, so we use the H\^opital's rule to get $\lim_{\alpha \to 0} ln  H_\alpha(\theta) = \lim_{\alpha \to 0} (\sum_i ln ((\epsilon(\theta) h(x_i|\theta))^{\lambda_i}) exp(\alpha ln (\epsilon(\theta) h(x_i|\theta)))/\sum_i \lambda_i   (\epsilon(\theta) h(x_i|\theta))^\alpha)
=\sum_i ln ((\epsilon(\theta) h(x_i|\theta)))^{\lambda_i}$. Note that, the H-C $lim_{\alpha \to 0}H_\alpha(\theta)$ is the geometric mean and it is the  Fisher's likelihood to the power $1/n$. It follows that the Fisher's likelihood inherits from the geometric mean the following feature:  it is dominated by smallest values of the set \{$\epsilon(\theta) h(x_1|\theta), \cdots, \epsilon(\theta) h(x_n|\theta)$\}. 

{\bf P2:} $H_\alpha (\theta)$ reduces to  the arithmetic mean when $\alpha=1$, and Harmonic mean when $\alpha=-1$. This is straightforward from Eq.~\ref{HCentrality}.

{\bf P3:} The function $H_\alpha(\theta)$ is increasing in $\alpha$. In the case of the H\"older mean, there are several proofs based on H\"older, Jensen, and Schwartz inequalities,   theory of convex functions, and critical point theory~\cite{Stolarsky96,Beckenbach50,Ku99,Micic12}.  They are too long and require advanced mathematical concepts.  We provide here a simplest proof using elementary calculus in the case of $H_\alpha(\theta)$ in Eq.~\ref{HCentrality}.   Let us consider the sign of the  derivative of $H_\alpha(\theta)$ in $\alpha$:
\begin{equation}
\frac{\partial H_{\alpha}(\theta)}{\partial \alpha} =  \frac{H_\alpha(\theta)}{\alpha} (\sum_i \lambda_i (\frac{h(x_i|\theta)}{H_\alpha(\theta)})^\alpha ln h(x_i|\theta) -  ln H_\alpha)
\end{equation}
Because $\sum_i \lambda_i (\frac{h(x_i|\theta)}{H_\alpha(\theta)})^\alpha=1$, so 
\begin{equation}
\frac{\partial H_{\alpha}(\theta)}{\partial \alpha} = \frac{H_\alpha(\theta)}{\alpha} \sum_i \lambda_i (\frac{h(x_i|\theta)}{H_\alpha(\theta)})^\alpha (ln h(x_i|\theta) - ln H_\alpha(\theta))
\end{equation}
Equivalently,
\begin{equation}
\frac{\partial H_{\alpha}(\theta)}{\partial \alpha} =  \frac{H_\alpha(\theta)}{\alpha} \sum_i \lambda_i (\frac{h(x_i|\theta)}{H_\alpha(\theta)})^\alpha ln \frac{h(x_i|\theta)}{H_\alpha(\theta)}
\end{equation}
Because the term $z_i = (h(x_i|\theta)/H_\alpha(\theta))^\alpha>0$, then $z_i>  ln(z_i)$ and therefore the lower band of the derivative  of $H_\alpha(\theta)$ is straightforward:
\begin{equation}
\frac{\partial H_{\alpha}(\theta)}{\partial \alpha} > H_\alpha(\theta) \sum_i \lambda_i  (ln \frac{h(x_i|\theta)}{H_\alpha(\theta)})^2 >0
\end{equation}
It follows that $H_\alpha(\theta)$ is an increasing function; that is  there is no critical point  $\alpha$ of $\frac{\partial H_{\alpha}(\theta)}{\partial \alpha}$. Consequently, the choice of $\alpha$ is an open issue that deserves further study.

{\bf P4:}  If the PDF $ h(x|\theta)$ is bounded  then $H_\alpha(\theta)$ is also a bounded function. Indeed,  $H_\alpha(\theta)$ is mean of $ \epsilon(\theta) h(x_1|\theta) \le  \cdots \le  \epsilon(\theta) h(x_n|\theta) $, so  $ \epsilon(\theta) h(x_1|\theta) \le H_\alpha(\theta)\le   \epsilon(\theta) h(x_n|\theta)$ by definition of the mean. Moreover, $lim_{\alpha \rightarrow + \infty} H_\alpha(\theta)= lim_{\alpha \rightarrow + \infty} \epsilon(\theta) h(x_n|\theta) (\sum_i \lambda_i (h(x_i|\theta) / h(x_n|\theta))^\alpha)^{1/\alpha}= \epsilon(\theta) h(x_n|\theta)$. Similarly, $lim_{\alpha \rightarrow - \infty} H_\alpha(\theta)= lim_{\alpha \rightarrow + \infty} \epsilon(\theta) h(x_1|\theta) (\sum_i \lambda_i (h(x_i|\theta) / h(x_1|\theta))^{-\alpha})^{1/-\alpha}= \epsilon(\theta) h(x_1|\theta)$.

{\bf P5:}  The H-C $H_{\alpha}(\theta)$  is a bijective function in $\alpha$.  Indeed, the inverse function theorem asserts that strictly increasing continuous functions are bijective and therefore invertible in the interval $[min_i  h(x_i|\theta), max_i  h(x_i|\theta)]$. 

{\bf P6:} The  data $X$ and therefore the PDF values $h(x_1|\theta), \cdots, h(x_n|\theta)$ do not have the same relevance in $H_\alpha (\theta)$. Indeed,  there are three data selection strategies  embedded in $H_\alpha (\theta)$. The first is the weight  $\lambda_i>0$  encoding the relative relevance of the $i^{th}$ data element; i.e $\sum_i \lambda_i=1$. The second is materialized by transform $\alpha$;  for a data element $x$,  its relevance $h^\alpha(x|\theta)$ is increasing or decreasing  depending $h(x|\theta)$ itself. The $1/\alpha$ transform does not cancel the $\alpha$ transform effects. The centrality $H_\alpha(\theta)$ 
is increasing or decreasing depending on $H^\alpha_\alpha(\theta)$. The third strategy concerns the influence of a PDF value $h(x_i|\theta)$ itself on the $H_\alpha(\theta)$. The reader can find more about the data selection strategies embedded in Pythagorean, H\"older, and Lehmer means in~\cite{Ziou23a,Ziou23b}.

%\subsection{Properties of  H-C as a probability}
%From a purely statistical point of view, the likelihood (i.e. $lim_{\alpha \to 0}H_{\alpha}^n(\theta)$) is not a probability~\cite{Fisher22}. This allows us to conclude that $H_{\alpha}^n(\theta)$ is not a probability. The H-C $H_{\alpha}(\theta)$ refers to an objective function, where its maximum is  the best statistical model for a given PDF and given data. It is a joint probability of the data $P_\theta(X)$, but it does not integrate to one with respect to $\theta$.   However, from the information theory, things can be seen differently.

{\bf P7:}  For a fixed $\theta$, the H-C $H_{\alpha}(\theta)$ is a probability of an observation falling within a cell $s \mp \epsilon(\theta)/2$. Indeed,  according to P3, $H_{\alpha}(\theta)  \in [min_i \epsilon(\theta) h(x_i|\theta), max_i \epsilon(\theta) h(x_i|\theta)]$.  Because the probability $\epsilon(\theta) h(x|\theta)$ is assumed to be continuous in $x$, it follows that there is at least one cell $s\mp\epsilon(\theta)/2$ of $x$ such that $H_{\alpha}(\theta) \simeq \epsilon(\theta) h(s|\theta) \in  [min_i \epsilon(\theta) h(x_i|\theta), max_i \epsilon(\theta) h(x_i|\theta)]$.

{\bf P8:}  For a fixed $\theta$, the cell $s \mp \epsilon(\theta)/2$ is unique if $h(x|\theta)$ is a bijective function in $x$.  Indeed, if  $h(x|\theta)$ is a bijective function in $x$, then the function inversion theorem ensures that there exists a unique $s=h^{-1}(H_{\alpha}(\theta)/\epsilon(\theta))$. 
If $h(x|\theta)$ is only surjective, there are $m \in ]1,n[$ cells  ${s_1, s_2, \cdots, s_m}$ such that $h(s_1|\theta)= ...= h(s_m|\theta)=H_{\alpha}(\theta)/\epsilon(\theta)$. It is worth noting that for a fixed $\theta$ and given data $X$, the parameter $\alpha$ defines completely the H-C and the cell $s$. To make the dependencies explicit, in the following, the center of the cell will be denoted $s_{\alpha,\theta}$.

{\bf P9:}  For a fixed $\theta$, let $\alpha_1 \le \alpha_2$, we have  $h(s_{\alpha_1,\theta}|\theta) \le h(s_{\alpha_2,\theta}|\theta)$. This is straightforward because according to {\bf P3}, we have  $H_{\alpha_1}(\theta) = \epsilon(\theta) h(s_{\alpha_1,\theta}) < H_{\alpha_1}(\theta) = \epsilon(\theta) h(s_{\alpha_2,\theta})$ and $\epsilon(\theta)>0$. It follows that $ h(s_{\alpha_1,\theta}) <  h(s_{\alpha_2,\theta})$.

{\bf P10:}  $\int_{-\infty}^{+\infty} H_{\alpha}(\theta) d \alpha =+\infty$. Indeed,  since $H_{\alpha}(\theta) \ge c=min_i \epsilon(\theta) h(x_i|\theta) > 0$, then $\int_{-\infty}^{+\infty} H_{\alpha}(\theta) d \alpha \ge c \int_{- \infty}^{+\infty}  d \alpha = +\infty $.

{\bf P11:} The expectation of $H_{\alpha}(\theta)$ with regards to a PDF $p(\alpha)$ is given by $\int_{-\infty}^{+\infty} p(\alpha) H_{\alpha}(\theta) d \alpha =  \epsilon(\theta) \\ \int_{-\infty}^{+\infty}  p(\alpha) h(s_{\alpha,\theta}|\theta)d \alpha$.   The equality is inferred from {\bf P8},  because a cell $s_{\alpha,\theta} \mp \epsilon(\theta)/2$ is chosen among $m \ge 1$ cells; i.e. $H_{\alpha}(\theta) =  h(s_{\alpha,\theta}|\theta)$.

%{\bf P11 : est-ce que le  $s_{\alpha,\theta}$ est une statistique suffisante? Reste a repondre a cette question.}

%%Indeed, because the logarithm is a concave function, thus  Jensen inequality leads to $ln \sum_i \lambda_i x^a_i  \ge  \sum_i\lambda_i ln x^a_i$.  The log H-likelihood has a lower band
%%\begin{equation} LH_\alpha (\theta) = \frac{1}{\alpha} ln \sum_i \lambda_i h(x_i|\theta)^\alpha \ge  \sum_i \lambda_i ln h(x_i|\theta)  \end{equation}
%The right term is log of the Fisher MLE. Let us now consider that $\theta^*$ is the MLE. At this point:
%\begin{equation} LH_\alpha (\theta) = \frac{1}{\alpha} ln \sum_i \lambda_i h(x_i|\theta^*)^\alpha \ge  \sum_i \lambda_i ln h(x_i|\theta^*)  \end{equation}
%This shows that the Fisher maximum likelihood  is not greater than the H-L.  

\section{Lehmer measure}
\label{LMeasureS}
From Eq.~\ref{Lmean}, we define the Lehmer centrality (L-C) as follow: 
\begin{equation}
 L_\alpha(\theta) = \epsilon(\theta) \frac{\sum_i \lambda_i  h^\alpha(x_i|\theta)}{\sum_i \lambda_i h^{\alpha-1}(x_i|\theta)}
 \label{LCentrality}
\end{equation}

As for  H-C, the L-C is a mean of the probabilities of an observation falling within a cell of size $\epsilon(\theta)$. The data selection mechanisms embedded in L-C are similar to those embedded in H-C (see {\bf P6}).  

{\bf Q1:} $L_0(\theta)$ and $ L_1(\theta)$ are related  the harmonic and arithmetic  mean of $h^\alpha(x_i|\theta),~i=1,n$ respectively. This is straightforward from Eq.~\ref{LCentrality}.

{\bf Q2:} The L-C is an increasing function in $\alpha$. Indeed,
\begin{equation}
 \frac{\partial L_{\alpha}(\theta)}{\partial \alpha}  = \epsilon(\theta) \frac{\sum_i \sum_j \lambda_i \lambda_j h^{\alpha-1}(x_i|\theta) h^{\alpha-1}(x_j|\theta) r_{ij}}{2 (\sum_i h^{\alpha-1}(x_i|\theta))^2} 
\end{equation}
where $ r_{ij}=( h(x_i|\theta) - h(x_j|\theta))(ln h(x_i|\theta) - ln h(x_j|\theta))$. Because the sign of   $h(x_i|\theta) -  h(x_j|\theta)$ and  ln $h(x_i|\theta) - ln h(x_j|\theta)$ is the same,  then $ r_{ij} \ge 0$.  All the other terms are positive and not all $r_{ij}$ are zeros. It follows that $\frac{\partial L_{\alpha}(\theta)}{\partial \alpha} >0$.

{\bf Q3:} {\bf P4} until {\bf P11} remains valid for the L-C.

\section{Comparison of H-C and L-C}
\label{CHL}

{\bf R1:}  H-C and L-C satisfy

\begin{equation}
 L_\alpha (\theta)\left \{
   \begin{array}{l c r}
       > H_\alpha(\theta)   & \mbox{if $\alpha >1 $} \\
      <   H_\alpha(\theta)  & \mbox{if $\alpha <1 $} \\
      = H_\alpha(\theta)  & \mbox{if $\alpha =1$} 
   \end{array}
   \right . 
   \label{Ordre}
\end{equation}
Indeed, the L-C can be rewritten as
\begin{equation}
 L_\alpha(\theta) = (\frac{H_\alpha(\theta)}{H_{\alpha-1}(\theta)})^{\alpha-1} H_\alpha(\theta) 
 \label{HLRelation}
\end{equation}
Because $H_\alpha(\theta)$ is increasing in $\alpha$, we conclude that  $(\frac{H_\alpha(\theta)}{H_{\alpha-1}(\theta)})^{\alpha-1}$ is greater (resp. smaller) than one if $\alpha \in ]1,+\infty[$ (resp. $\alpha \in ]-\infty,1[$). Substituting $\alpha=1$ in  Eqs.~\ref{HCentrality} and \ref{LCentrality} leads to $H_1(\theta)= L_1 (\theta)$. From {\bf P4} and {\bf Q3}, $H_\alpha(\theta)=L_{\alpha}(\theta)=max_i \epsilon(\theta) h(x_i|\theta)$ (resp. $H_\alpha(\theta)=L_{\alpha}(\theta)=min_i \epsilon(\theta) h(x_i|\theta))$ for $\alpha \to \infty$ (resp. $\alpha \to -\infty$).

\section{Maximum centrality estimation}
\label{MCE}
Having defined the H\"older and Lehmer centralities $H_\alpha(\theta)$ and $L_\alpha(\theta)$, we will now deal with the estimation of $\theta$, by maximizing each of they. None of these measures has a maximum in $\alpha$ because they are increasing with respect to this variable. Consequently, for now we assume that $\alpha$ is known. Moreover, let us recall  that  $H_\alpha(\theta)$ and $L_\alpha(\theta)$ are assumed both derivable in $\theta$. In what follows, for the sake of simplicity, the cell size $\epsilon(\theta)$ is ignored and the logarithm of both measures are designed by $C_\alpha(\theta)$ where $C=\ln(H)$ in the case of H\"older centrality  and $C=\ln(L)$ for Lehmer one.
  
{\bf C1:}  $C_\alpha(\theta)$ has maximum in $\theta$ if $ h(s_{\alpha,\theta}|\theta)$ has a maximum in $\theta$. Indeed, $C_\alpha(\theta)= ln(h(s_{\alpha,\theta}|\theta))$ and the PDF is supposed to have a maximum  with respect to $\theta$. Note that $s_{\alpha,\theta}$ depends on the value of $C$ being either $H$ or $L$.

%{\bf C2:} Neither H-C nor L-C has a  critical point on $\alpha$. This is straightforward from P1 and Q1. However, according to P9 and Q9, a judicious choice of a prior $p(\alpha)$  could lead to maximizing H-C or L-C.

{\bf C2:} Let us now, look for  $\theta$ as a critical point of the C-centrality.
\begin{equation}
 C_\alpha (\theta) = \left \{
   \begin{array}{l r }
      \ln H_\alpha^\alpha(\theta) -  ln H_{\alpha-1}^{\alpha-1}(\theta)  & \mbox{if $C=L  $} \\
      \frac{1}{\alpha} \ln H_\alpha^\alpha(\theta)   & \mbox{if $C=H  $ }     
   \end{array}
   \right . 
\end{equation}
Note that, $ H_\alpha^\alpha(\theta)|_{\alpha=0}=1$. The first order condition w.r.t $\theta$ is:
\begin{equation}
\frac{\partial { C_\alpha}(\theta)}{\partial \theta} =
\left \{
   \begin{array}{l r }
       \frac{1}{H_\alpha^\alpha} \frac{\partial H_\alpha^\alpha(\theta)}{\partial \theta} -   \frac{1}{H_{\alpha-1}^{\alpha-1}(\theta)} \frac{\partial H_{\alpha-1}^{\alpha-1}(\theta)}{\partial \theta}=0  & \mbox{if $C=L  $} \\
      \frac{1}{\alpha H_\alpha^\alpha(\theta)} \frac{\partial H_\alpha^\alpha(\theta)}{\partial \theta} =0   & \mbox{if $C=H  $} 
   \end{array}
   \right . 
\label{FirstOrderC}
\end{equation}
Substituting H-L in Eq.~\ref{HCentrality}, we obtain:
\begin{equation}
\frac{\partial { C_\alpha}(\theta)}{\partial \theta} =
\left \{
   \begin{array}{l r }
      \alpha E_{g_\alpha}(\frac{\partial ln h(x|\theta)}{\partial \theta}) -  (\alpha-1) E_{g_{\alpha-1}}(\frac{\partial ln h(x|\theta)}{\partial \theta})=0  & \mbox{if $C=L  $} \\
      E_{g_\alpha}(\frac{\partial ln h(x|\theta)}{\partial \theta}) =0   & \mbox{if $C=H  $} 
   \end{array}
   \right . 
\label{FirstOrderC1}
\end{equation}
where  $E_{g_\alpha}(\frac{\partial ln h(x|\theta)}{\partial \theta})=\sum_i \frac{\partial ln h(x_i|\theta)}{\partial \theta} g_\alpha(x_i|\theta)$ and the PDF $g_\alpha(x_i|\theta)= \lambda_i h^\alpha(x_i|\theta)/H_\alpha^\alpha(\theta)$.

{\bf C4:} The critical points $\theta^*_L$ and $\theta^*_H$ of H-C and L-C fulfills  $\theta^*_L=\theta^*_H$ if and only if $\alpha=1$ or when $\theta^*_H$ is a critical point of both $H_\alpha^\alpha(\theta)$ and   $H_{\alpha-1}^{\alpha-1}(\theta)$. This is straightforward from  Eq.~\ref{FirstOrderC1}.

{\bf C5:} The second order condition w.r.t $\theta$ is given by:
\begin{equation}
\frac{\partial^2 C_\alpha (\theta)}{\partial \theta^2}  =  
\left \{
   \begin{array}{l r }
      \frac{1 }{H_\alpha^\alpha(\theta)} \frac{\partial^2 H_\alpha^\alpha(\theta)}{\partial^2 \theta}- (\frac{1}{H_\alpha^\alpha(\theta)} \frac{\partial H_\alpha^\alpha(\theta)}{\partial \theta})^2-  \frac{1}{H_{\alpha-1}^{\alpha-1}(\theta)} \frac{\partial^2 H_{\alpha-1}^{\alpha-1}(\theta)}{\partial^2 \theta}+ (\frac{1}{H_{\alpha-1}^{\alpha-1}(\theta)} \frac{\partial H_{\alpha-1}^{\alpha-1}(\theta)}{\partial \theta})^2 & \mbox{if $C=L  $} \\
       \frac{1 }{\alpha} \frac{1 }{H_\alpha^\alpha(\theta)} \frac{\partial^2 H_\alpha^\alpha(\theta)}{\partial^2 \theta}- \alpha (\frac{1}{\alpha H_\alpha^\alpha(\theta)} \frac{\partial H_\alpha^\alpha(\theta)}{\partial \theta})^2 & \mbox{if $C=H  $} 
   \end{array}
   \right . 
\end{equation}
At the critical point $\theta$ (Eq.~\ref{FirstOrderC}), the second derivative is rewritten:
\begin{equation}
\frac{\partial^2 C_\alpha (\theta)}{\partial \theta^2}  =   \left \{
   \begin{array}{l r }
      \frac{1 }{H_\alpha^\alpha(\theta)} \frac{\partial^2 H_\alpha^\alpha(\theta)}{\partial^2 \theta}-  \frac{1}{H_{\alpha-1}^{\alpha-1}(\theta)} \frac{\partial^2 H_{\alpha-1}^{\alpha-1}(\theta)}{\partial^2 \theta}  & \mbox{if $C=L  $} \\
    \frac{1 }{\alpha H_\alpha^\alpha(\theta)} \frac{\partial^2 H_\alpha^\alpha(\theta)}{\partial^2 \theta}  & \mbox{if $C=H  $} 
   \end{array}
   \right . 
\label{SecondOrderC}
\end{equation}
We have 
\begin{equation}
  \frac{\partial^2 ln H_\alpha^\alpha(\theta)}{\partial \theta^2}  =  \frac{\frac{\partial^2  H_\alpha^\alpha(\theta) }{\partial \theta^2}}{H_\alpha^\alpha(\theta)}-(\frac{\partial ln H_\alpha^\alpha(\theta)}{\partial \theta})^2
\end{equation}
Substituting this equation in Eq.~\ref{SecondOrderC} and by using $ \frac{\partial ln H_\alpha^\alpha(\theta)}{\partial \theta} =  \frac{\partial ln H_{\alpha-1}^{\alpha-1}(\theta)}{\partial \theta}$ when C=L and $\frac{\partial ln H_\alpha^\alpha(\theta)}{\partial \theta}=0$ when C=H, the second derivative becomes:
\begin{equation}
\frac{\partial^2 C_\alpha (\theta)}{\partial \theta^2}  =   \left \{
   \begin{array}{l r }
      \frac{\partial^2 ln H_\alpha^\alpha(\theta)}{\partial \theta^2}-   \frac{\partial^2 ln H_{\alpha-1}^{\alpha-1}(\theta)}{\partial \theta^2}  & \mbox{if $C=L  $} \\
      \frac{1 }{\alpha} \frac{\partial^2 ln H_\alpha^\alpha(\theta)}{\partial \theta^2}  & \mbox{if $C=H  $} 
   \end{array}
   \right . 
\label{2Deriv}
\end{equation}
Straightforward manipulations leads to
\begin{equation}
\frac{\partial^2 C_\alpha (\theta)}{\partial \theta^2}  =   \left \{
   \begin{array}{l r }
      \alpha E_{g_\alpha} ((\frac{\partial ln h(x|\theta)}{\partial  \theta})^2) +  E_{g_\alpha}  (\frac{\partial^2 ln h(x|\theta)}{\partial  \theta^2}) - (\alpha-1) E_{g_{\alpha-1}} ((\frac{\partial ln h(x|\theta)}{\partial  \theta})^2) -  E_{g_{\alpha-1}}  (\frac{\partial^2 ln h(x|\theta)}{\partial  \theta^2})   & \mbox{if $C=L  $} \\
        \alpha E_{g_\alpha} ((\frac{\partial ln h(x|\theta)}{\partial  \theta})^2) +  E_{g_\alpha}  (\frac{\partial^2 ln h(x|\theta)}{\partial  \theta^2})  & \mbox{if $C=H  $} 
   \end{array}
   \right . 
\label{2DerivDetails}
\end{equation}
In the case of H-C, the critical point is a maximum if  $E_{g_\alpha}  (\frac{\partial^2 ln h(x|\theta)}{\partial  \theta^2}) < - \alpha E_{g_\alpha} ((\frac{\partial ln h(x|\theta)}{\partial  \theta})^2) < 0$ when $\alpha \ge 0$ and  $E_{g_\alpha}  (\frac{\partial^2 ln h(x|\theta)}{\partial  \theta^2}) < - \alpha E_{g_\alpha} ((\frac{\partial ln h(x|\theta)}{\partial  \theta})^2)$ when $\alpha < 0$. In the case of L-C, the critical point is a maximum if  $\alpha (E_{g_\alpha} ((\frac{\partial ln h(x|\theta)}{\partial  \theta})^2) - E_{g_{\alpha-1}} (\frac{\partial ln h(x|\theta)}{\partial  \theta})^2))  
   <    E_{g_{\alpha-1}}  (\frac{\partial^2 ln h(x|\theta)}{\partial  \theta^2}) - E_{g_\alpha}  (\frac{\partial^2 ln h(x|\theta)}{\partial  \theta^2}) - E_{g_{\alpha-1}} ((\frac{\partial ln h(x|\theta)}{\partial  \theta})^2)$.

%{\bf C6:}  The values of $H_\alpha(\theta^*_H)$ and $L_\alpha(\theta^*_L)$, if any, increase with $\alpha$ without exceeding $max_i \epsilon(\theta) h( x_i |\theta)$. This is straightforward from P1 and Q1.

\section{Accuracy of C-estimators}
\label{AE}
The maximum likelihood estimators are often studied according to Fisher's requirements which are consistency, efficiency and sufficiency. But instead of the theoretical study, we propose to define easy-to-implement performance measures for the evaluation of the accuracy of the distribution fit. Even several are possible to define, in this study, we will limit ourselves to two measures.  For now, we assume that the parameter $\alpha$ is known.

\subsection{Residuals}
 For a given observation $x_i$ and a model $\theta$, the residual can be measured by a norm of the absolute value of the difference between the observation frequency $d_i$ and its predicted value $h(x_i|\theta)$. For all observations, the arithmetic mean of the residuals is often used~\cite{NIST12} for the evaluation of the C-estimator $\theta$. Rather, we propose the use of H\"older and Lehmer mean families instead. Let's start with H\"older  with the parameter $\beta  \in \mathbb{R}$:
\begin{equation}
e_H = (\sum_{i=1}^n  \lambda_i  |d_i  - h(x_i|\theta)|^\beta)^{1/\beta}
\label{ResidualErrorH}
\end{equation} 
When $\beta=2$, $ e_H$ reduces to the famous square root of the mean square errors. The residual error $e_H$ inherits the properties of the H\"older centrality measure: it is increasing in $\beta$ and implements a data selection mechanisms~\cite{Ziou23b,Ziou23a}.
Let us now approach the interpretation of $e_H$ in terms of probability. Assuming that the error $\eta = \lambda (d  - h(x|\theta))$ can be sampled from a normal distribution $N(0,\sigma)$, then the pdf of $|\eta|$  is a half-normal with mean $\sigma \sqrt{2/\pi}$ and a variance $\sigma^2 (1-2/\pi)$.  The pdf of $|\eta|^\beta$ is straightforward; let $Y=|\eta|^{\beta}$,  $Pr(Y \le y)=Pr(|\eta|^{\beta} \le y) = Pr(|\eta| \le y^{1/\beta})=F_{|\eta|}(y^{1/\beta})$, where $F$ is the cumulative distribution function of the half-normal. 
Now, let us consider  the iid random variables $Y_1, \cdots, Y_n$, then according to the central limit theorem, the pdf of  $Z=\sum_i Y_i$ is normal distribution with the mean equals to $n E(Y)$ and variance $n Var(Y)$. Because the variable $Z$ is non negative, we then use the truncated normal distribution. Under this choice,  the pdf of $e_H$ is straightforward; let $\psi =Z^{1/\beta}$,  $Pr(\psi \le z)=Pr(Z^{1/\beta} \le z) = Pr(Z \le z^\beta)=F_{Z}(z^\beta)$, where:
\begin{equation}
	f_H(z)= F_Z^{'}(z^{\beta})=\beta z^{\beta-1} g(z^{\beta}|n E(Y),n Var(Y) ) / (1/2+erf ( E(Y)/ Var(Y)))
\end{equation}
and $g$ is the normal PDF and   $z \ge 0$. 

For Lehmer   with the parameter $\beta \in \mathbb{R}$, we define the residual  error as follow: 
\begin{equation}
e_L= \frac{\sum_i  \lambda_i |d_i  - h(x_i|\theta)|^{\beta}}{\sum_i  \lambda_i  |d_i  - h(x_i|\theta)|^{\beta-1}} 
\label{ResidualErrorL}
\end{equation}
The residual error  $e_L$ inherits the properties of the Lehmer centrality measure: it is increasing in $\beta$ and implements data selection mechanisms~\cite{Ziou23b,Ziou23a}.  The random variable $Z$ associated with the error measure $e_L$ is the ratio of two dependent truncated normal random variables. We propose the use of an easy-to-implement approximation by considering that both variables are normal and independent having means $\mu_\beta$ and $\mu_{\beta-1}$ and standard deviations $\sigma_\beta$ and $\sigma_{\beta-1}$. In this case, the PDF of the ratio of two normal distributions is given in~\cite{Diaz-Frances13}. Because the variable $Z$ is non negative, we then use its truncated version given by:
\begin{equation}
	f_L(z)=C(\rho, \mu,\delta_{\beta-1}) \frac{\rho}{\pi(1+\rho^2 z^2)} exp(-\frac{\rho^2 \mu^2 +1}{2 \delta_{\beta-1}^2} ) (1 + \sqrt{\frac{\pi}{2}}q erf(\frac{q}{\sqrt{2}}) exp(\frac{q^2}{2}) )
\end{equation}
where $z \ge 0$, $C(\rho, \mu,\delta_{\beta-1})$ a normalization constant, $q=(1+\mu \rho^2 z)/\delta_{\beta-1} \sqrt{1+\rho^2 z^2}$, $\rho=\sigma_{\beta-1}/\sigma_{\beta-1}$, $\mu=\mu_{\beta}/\mu_{\beta-1}$, and  $\delta_{\beta-1}=\sigma_{\beta-1}/\mu_{\beta-1}$.

%Further we will study, $H_R = (\sum_i  \lambda_i  |d_i  - h(x_i|\theta)|^\beta )^{1/\beta}$.

\subsection{Observed shape} 
We propose to define a uncertainty measure of the C-estimator $\theta^*$ inspired by Fisher information. For a PDF $h(x|\theta)$, the Fisher information about  $\theta$ is  the expected value of second derivative of $-ln(h(x|\theta))$ with regards to $\theta$, evaluated at  $\theta^*$.  The interpretation is as follows.  The change of the function $ln(h(x|\theta))$ on the parameter $\theta$ space in the vicinity of $\theta^*$ provides information about the uncertainty of the C-estimator  $\theta^*$. For example, the change in orientation of the tangent of $ln(h(x|\theta))$ in the vicinity of $\theta^*$ provides information about the shape of the peak of $ln(h(x|\ theta))$ at the point  $\theta^*$. If the change is rapid, then for any $\theta$ belonging to the vicinity of $\theta^*$, $ln(h(x|\theta^*))$ can be easily distinguished from any $ln( h(x|\theta))$. In this case, the accuracy of the estimated  $\theta^*$, the second derivative of $ln(h(x|\theta))$ with regards to $\theta$ at the point $\theta^*$, and therefore the Fisher information are all high. The second derivative of $ln(h(x|\theta))$ with regards to $\theta$ at the point $\theta^*$ is the curvature measure at this point. This geometric interpretation highlights the strength of Fisher's information which establishes an equivalence between uncertainty and this geometric characteristic.   It happens that the arithmetic mean is used instead of the expectation and the Fisher information, known as observed Fisher information, is the  negative of the second derivative  of the logarithm of the likelihood. For the estimation of centrality, we take inspiration from the observed Fisher information to propose observed C-Fisher information. Since the uncertainty of the C-estimator $\theta*$ depends on the curvature of the measured centrality, we define observed C-Fisher information, as the second derivative of the logarithm of minus the centrality given in Eq.~\ref{2DerivDetails}.  Equivalently, because the centrality is equal to $ln h(s_{\alpha,\theta}|\theta)$ (See {\bf P7} and {\bf Q7}), then  the observed C-Fisher information, can be seen also as the second derivative of the  $-ln h(s_{\alpha,\theta}|\theta)$.

%\subsection{Estimation of $\alpha$} Both H-C and L-C are increasing wrt the parameter $\alpha$. It can be estimated by maximizing  H-C or L-C only if constraints are added delimiting the useful values.  Following {\bf P11} and {\bf Q11}, it can be estimated by using specific prior knowledge.  The prior $p(\alpha)$ indicates the more likely values of $\alpha$. Another strategy consists in choosing $\alpha$ leading to the smallest the residual errors in~Eqs~\ref{ResidualErrorH} and~\ref{ResidualErrorL}.

\section{Case study}
\label{CS}

\begin{figure}[ht]
	%\centering
	\includegraphics[scale=0.45]{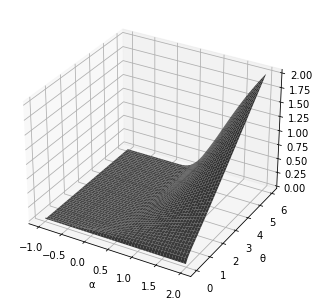}{\hspace*{-5mm} a­} 
	\includegraphics[scale=0.3]{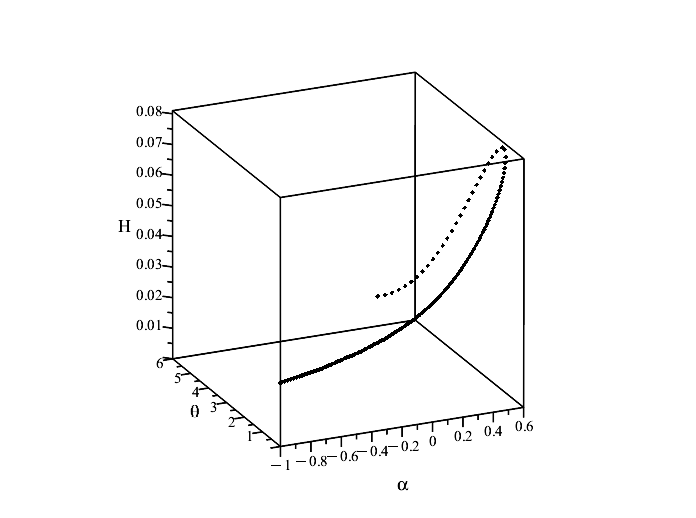}{\hspace*{-10mm} b} 
	\includegraphics[scale=0.3]{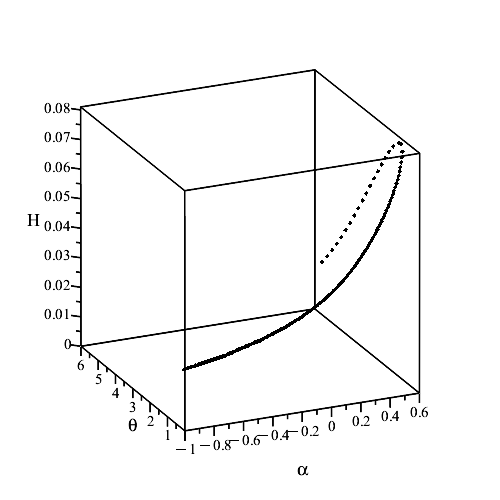}{\hspace*{-5mm} c} 
	\includegraphics[scale=0.35]{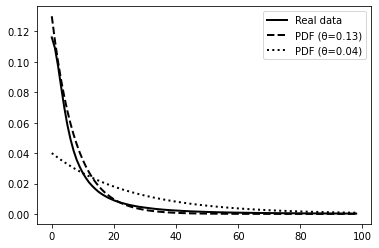}{d}
	\includegraphics[scale=0.45]{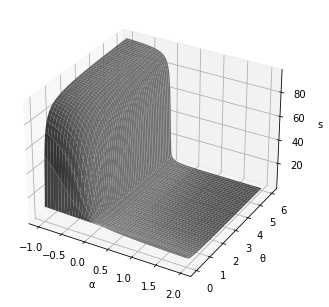}{e}  
	\includegraphics[scale=0.25]{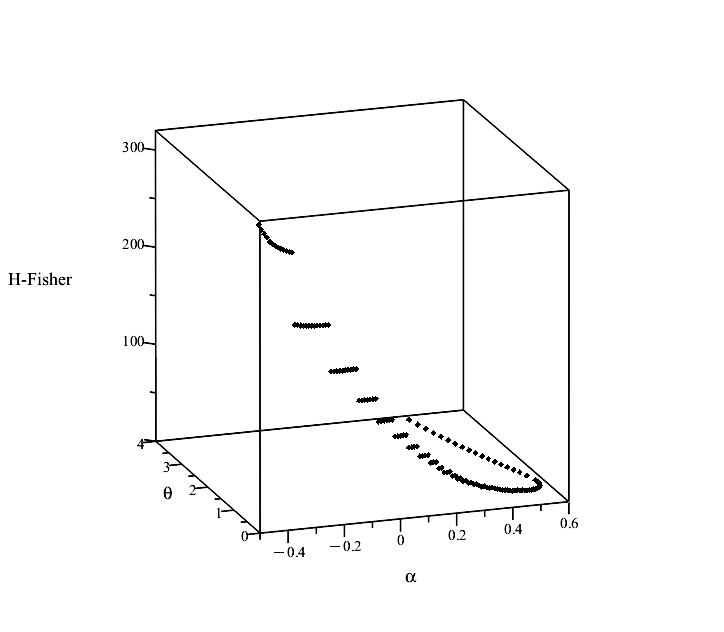}{f} 
	\caption{First row:  The H-C  calculated using the histogram (d),  critical points, and  maximums.  Second row:  a histogram, central observations explained in {\bf P7-P9}, and H-Fisher at the maximums.}
	\label{Critic}
\end{figure}

\begin{figure}[ht]
	%\centering
	\includegraphics[scale=0.45]{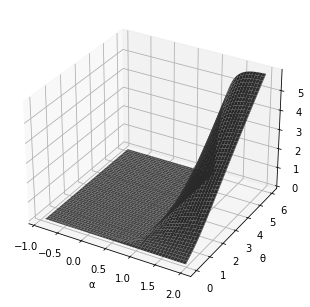}{a} 
	\includegraphics[scale=0.3]{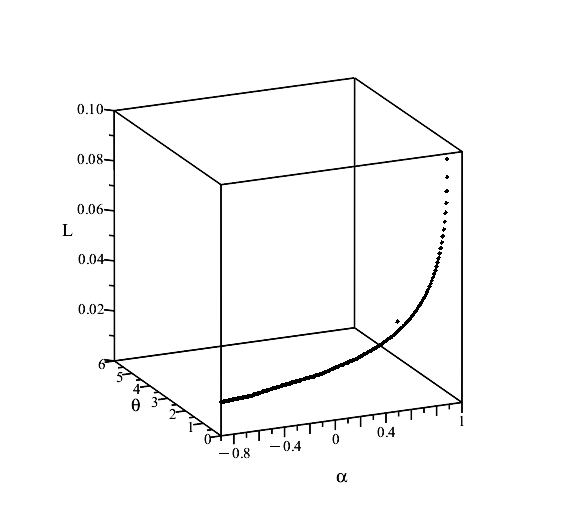}{b} 
	\includegraphics[scale=0.3]{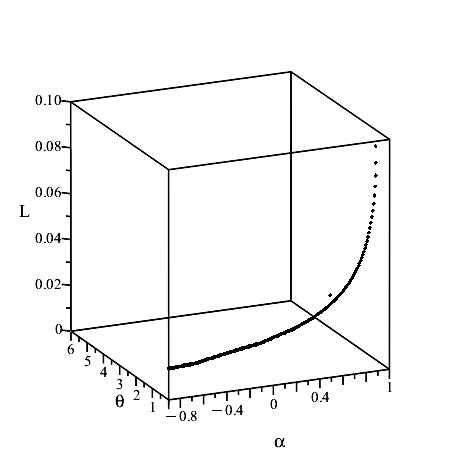}{c} 
	\centering
	\includegraphics[scale=0.35]{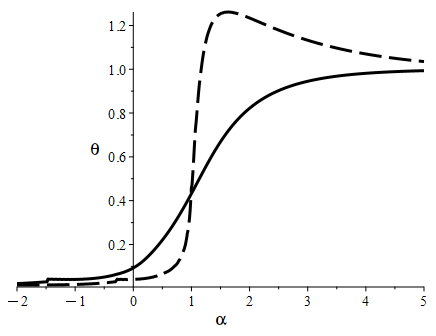}{h}   
	\caption{First row:   The L-C   calculated using the histogram in~\ref{Critic}.d.,  critical points, and  maximums. Second row:  the critical points  $\theta$ of H\"older (solid) and  Lehmer (dash) centralities as function of $\alpha$ drawn using the implicit function in Eq.~\ref{ExpFirstCond}.}
	\label{CriticL}
\end{figure}

We propose to study the C-estimator of the exponential PDF $h(x|\theta) = \theta exp(-\theta x)$. This bijective PDF in $x$ was widely used in several areas including  image processing, queuing theory, and physics~\cite{Kerouh18,Sundarapandian9,Berger19}.  The first derivative of centrality in Eq.~\ref{FirstOrderC} is rewritten as follows:
\begin{equation}
\frac{\partial C_\alpha(\theta)}{\partial \theta} =    \left \{
   \begin{array}{l r }
    \sum_i \lambda_i (- x_i + \frac{1}{\theta}) (\alpha  g_\alpha(x_i|\theta)  -  (\alpha-1)  g_{\alpha-1}(x_i|\theta))  & \mbox{if $C=L  $} \\
        \sum_i \lambda_i (- x_i + \frac{1}{\theta})  g_\alpha(x_i|\theta)  & \mbox{if $C=H  $} 
   \end{array}
   \right . 
\label{ExpFirstCond0}   
\end{equation}
where $g_\alpha(x|\theta)= exp(-\alpha \theta x)/\sum_i \lambda_i exp(-\alpha \theta x_i)$. Hence, the critical point is the following  fixed point: 
\begin{equation}
\theta =  \left \{ \begin{array}{l r }
\frac{ 1}{ \sum_i \lambda_i x_i  (\alpha  g_\alpha(x_i|\theta)  -  (\alpha-1)  g_{\alpha-1}(x_i|\theta))} & \mbox{if $C=L  $} \\
\frac{ 1 }{ \sum_i  \lambda_i x_i  g_\alpha(x_i|\theta)}   & \mbox{if $C=H  $}
 \end{array}
   \right . 
\label{ExpFirstCond}
\end{equation}
Fig.~\ref{Critic}.d depicts a histogram  of the absolute value of the DCT coefficients extracted from a gray-level image. The exponential of the centrality H-C, its critical points, and its maximums are shown in Figs.~\ref{Critic}.a,~\ref{Critic}.b, and~\ref{Critic}.c as function of $\alpha$ and $\theta$. Tow observations can be made. First, these functions are increasing in $\alpha$ as it is expected ({\bf P3}).  Second, some of critical points are not maximums as it is explicit in Figs.~\ref{Critic}.b and~\ref{Critic}.c. Fig.~\ref{Critic}.e depicts the central observation, denoted $s$ ({\bf P7-P9}), as function of $\alpha$ and $\theta$. As expected, a value of $s$ is lower for higher $\alpha$, because as $\alpha$ increases, the centrality increases ({\bf P3}) and hence in the case the centrality $h(s| \theta)$ ({\bf P7}), which is the exponential PDF, is higher at a lower $s$. The exponential of L-C in Eq.~\ref{LCentrality}, the critical points, and the maximums are plotted in Figs.~\ref{CriticL}.a,~\ref{CriticL}.b, and~\ref{CriticL}.c.
 Fig.~\ref{CriticL}.h shows the plots of the critical points $\theta$  in Eq.~\ref{ExpFirstCond} as functions of $\alpha$, using the histogram in Fig.~\ref{Critic}.d.  It can be noted that, the change in the value of $\theta$ caused by $\alpha$ is faster in the Lehmer case.  Fig.~\ref{CriticL}.h shows that there are estimators that can be derived from L-C but not from H-C. Moreover, the comparison of the two C-estimators in Eq.~\ref{ExpFirstCond} leads to conclude that for $\alpha>1$ (resp. $\alpha<1$) the Lehmer estimator (C=L) is higher (resp. lower) than that of H\"older (C=H).
 
We will now deal with second derivative of the two centralities at the critical points. Let us set $v(\theta,\alpha)=\sum_i \lambda_i x_i^2 g_\alpha(x_i|\theta) - (\sum_i \lambda_i x_i g_\alpha(x_i|\theta))^2$, according to Eq .~\ref{2DerivDetails}, the second derivatives at the critical point are given by:
\begin{equation}
\frac{\partial^2 C_\alpha}{\partial \theta^2} =  \left \{ \begin{array}{l r } \frac{-2\alpha}{\theta^2}+ \alpha^2  \sum_i \lambda_i x_i^2 g_\alpha(x_i|\theta) -  (\alpha-1)^2  \sum_i \lambda_i x_i^2 g_{\alpha-1}(x_i|\theta) & \mbox{if $C=L  $} \\
 - \frac{1 }{\theta^2}+ \alpha   v(\theta,\alpha) & \mbox{if $C=H   $}
\end{array}
   \right . 
 \label{ExpSecondCond}
\end{equation}
Note that,  $v(\theta,\alpha)> 0$. Indeed, let consider the positive function $\sum_i \lambda_i g_\alpha(x_i|\theta) (x_i - \sum_j \lambda_j x_j g_\alpha(x_j|\theta))^2$. It can be written as $\sum_i \lambda_i g_\alpha(x_i|\theta) x_i^2 - 2 \sum_i \lambda_i g_\alpha(x_i|\theta) x_i \sum_j \lambda_j x_j g_\alpha(x_j|\theta)+(\sum_j \lambda_j x_j g_\alpha(x_j|\theta))^2$. Simplifying it gives
 $\sum_i \lambda_i  x_i^2 g_\alpha(x_i|\theta) - (\sum_i \lambda_i x_i g_\alpha(x_i|\theta))^2=v(\theta,\alpha)$. 
 When it exists, the critical point of the H-C is a maximum if and only if   $ \alpha v(\theta,\alpha) < \frac{1 }{\theta^2}$. More specifically,  there is always a maximum when $\alpha \le 0$ or when  $\alpha v(\theta,\alpha)  \in [0, \frac{1 }{\theta^2}[$. Following the same reasoning for the case of L-C; the critical point is a maximum if 
$   \sum_i \lambda_i x_i^2 (\alpha^2g_\alpha(x_i|\theta) -   (\alpha-1)^2 g_{\alpha-1}(x_i|\theta)) < \frac{2\alpha}{\theta^2}$. From this inequality, we  deduce that if the left term is positive then $\alpha$ is positive and when $\alpha < 0$ then the left term is negative. About the uncertainty, the H-Fisher information estimated at the maximums is plotted Fig.~\ref{Critic}.f as function of $\alpha$ and $\theta$. The uncertainty is lower for higher $\alpha$. Another important question concerns the comparison of the C-estimator and the maximum likelihood estimator. We answer this question through an example in the case of H-Estimator. By using  the histogram in Fig.~\ref{Critic}.d, the best H-estimator according to the residuals $e_H$ in Eq.~\ref{ResidualErrorH} when $\beta=2$  is $\theta=0.1300$ obtained at $\alpha=0.0800$ and leading to $e_H=0.0005$. For comparison purposes, the maximum likelihood estimator (i.e. $\alpha=0$) is $\theta=0.04$, which gives $e_H=0.0214$ for $\beta=2$. The estimated PDFs are plotted in  Fig.~\ref{Critic}.d.

%\begin{table}
%\begin{tabular}{|l|c|c|c|c|c|} \hline
%$\alpha$ & ln H-L & $\hat{\theta}$ &  Shape & Residu \\ \hline  \hline
%-0.5 & -3.6239 & 0.0600 &  469.9815 & 0.0055 \\ \hline
%0 & -3.4107 & 0.09 &  123.4568 & 0.0114 \\ \hline
%0.5 & -3.1690 & 0.1500 &  25.7246 &  0.0337 \\ \hline
%\end{tabular}
%\caption{???}
%\label{PerfCasReal}
%\end{table}

\begin{figure}[ht]
\centering
\includegraphics[scale=0.35]{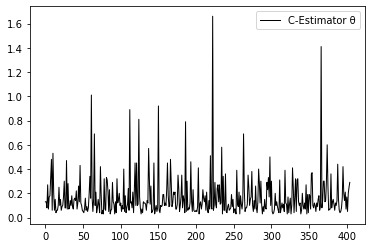}.a 
\includegraphics[scale=0.35]{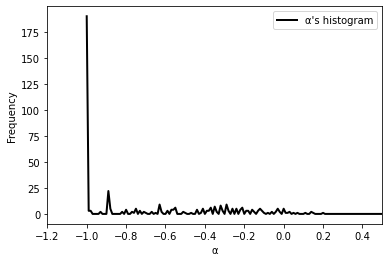}.b 
\caption{Fitting of 405 histograms of DC coefficients. a) The H-estimator $\theta$ for each image. b) The  histogram of the $\alpha$ selected according to the residual error.  }
\label{405}
\end{figure}

We end this case study by presenting an example of histogram adjustment of the absolute values of the DCT coefficients. To do this, we examine the distribution of the parameter $\alpha$ leading to the best H-estimator according to $e_H()$ in Eq.~\ref{ResidualErrorH}. Fig.~\ref{405} presents the best H-estimators of 405 histograms of DCT coefficients extracted from gray-level images. In Fig.~\ref{405}.a, an abscissa refers to a histogram and the corresponding ordinate is the value of the best H-estimator. For a histogram, the maximums of H-C are located and the one with the smallest $e_H()$ is chosen. The average $e_H()$ for the 405 histograms is 0.0008. The histogram of the variable $\alpha$ used to calculate the 405 best H-estimators is shown in~Fig.~\ref{405}.b. It can be observed that most of the values of $\alpha$ are not zero. Which indicates that, according to $e_H()$, maximum likelihood is not the best estimator.

\bibliographystyle{elsarticle-num}

%\section*{References}

\end{document}